\documentclass{amsart}
\newtheorem{theorem}{Theorem}
\newtheorem{proposition}{Proposition}
\newtheorem{lemma}{Lemma}
\newtheorem{corollary}{Corollary}
\title{Heegaard gradient of Seifert fibered 3-manifolds} 
\author{Kazuhiro Ichihara}
\address{Department of Information and Computer Sciences, 
Nara Women's University, Kita-Uoya Nishimachi, Nara 630-8506}
\email{ichihara@vivaldi.ics.nara-wu.ac.jp}
\thanks{The author is partially supported by JSPS Research Fellowships for Young Scientists.}
\begin{document}
\maketitle
\begin{abstract}
The infimal Heegaard gradient of a 3-manifold 
was defined and studied by Marc Lackenby 
in an approach toward the well-known virtually Haken conjecture. 
As instructive examples, we consider Seifert fibered 3-manifolds, and 
show that a compact orientable Seifert fibered 3-manifold has 
zero infimal Heegaard gradient if and only if 
it virtually fibers over the circle or 
over a surface other than the 2-sphere, 
equivalently, it has infinite fundamental group. 
\end{abstract}
%
\section{Introduction}\label{sec-intro}

The virtually Haken conjecture \cite[Problem 3.2]{K} 
is one of the most important conjectures 
in the study of 3-manifolds. 
Recently, in \cite{L}, 
Lackenby proved the following toward the virtually Haken conjecture. 
Let $M$ be a compact orientable 3-manifold 
with boundary a (possibly empty) union of tori. 
If the fundamental group of $M$ 
fails to have Property ($\tau$) with respect to 
some lattice of finite regular covers of $M$ and 
the lattice of covers has non-zero 
infimal (strong) Heegaard gradient, 
then $M$ is virtually Haken. 

The infimal Heegaard gradient 
was introduced by Lackenby in \cite{L} as follows. 
Recall that a \textit{Heegaard surface} 
in a compact orientable 3-manifold means 
a closed embedded surface 
which separates the manifold into two compression bodies. 
It is well-known that every compact orientable 3-manifold 
contains a Heegaard surface. 
Let $\{ M_i \to M \}$ be 
a collection of finite coverings of 
a compact orientable 3-manifold $M$ with covering degree $d_i$. 
Then 
the \textit{infimal Heegaard gradient of the collection $\{ M_i \to M \}$} 
is defined as 
$$ \inf_i \frac{\chi^h_- (M_i)}{d_i} \ ,$$
where $\chi^h_- (M_i)$ denotes 
the negative of the maximal Euler characteristic 
of a Heegaard surface in $M_i$. 
The \textit{infimal Heegaard gradient} of $M$ is defined as 
the infimal Heegaard gradient of 
the collection of all finite coverings of $M$. 
For brevity, we sometimes drop the word `infimal'. 

Toward the virtually Haken conjecture, 
on account of the Lackenby's result, 
our attention is attracted 
to the vanishing of Heegaard gradient. 
There actually exist 3-manifolds with 
zero Heegaard gradient. 
Easy examples are given by 
surface bundles over the circle $S^1$, 
all of which are easily shown 
to have zero Heegaard gradient. 
In fact, Lackenby proposed a conjecture: 
a compact orientable hyperbolic 3-manifold has 
zero infimal Heegaard gradient 
if and only if it virtually fibers over the circle 
(Heegaard gradient conjecture \cite{L}). 
He showed in \cite{L} this is also the case 
for closed orientable \textit{reducible} 3-manifolds, 
i.e., 3-manifolds which contain an embedded 2-sphere 
not bounding a 3-ball. 

In order to analyze when Heegaard gradient vanishes, 
we consider Seifert fibered 3-manifolds as instructive examples 
and study which (covers of) Seifert fibered 3-manifolds have 
zero Heegaard gradient. 
A \textit{Seifert fibered 3-manifold} is defined as 
a 3-manifold with a decomposition into disjoint circles, 
called \textit{fibers}, such that 
each circle has a neighborhood isomorphic to 
a fibered solid torus or Klein bottle. 
See \cite{S} for the basic theory of Seifert fibered 3-manifolds. 
Then our main result is: 
\begin{theorem}\label{thm-SF}
A compact orientable Seifert fibered $3$-manifold 
has zero infimal Heegaard gradient if and only if 
it virtually fibers over the circle or 
over a surface other than the $2$-sphere. 
Equivalently, 
a compact orientable Seifert fibered $3$-manifold 
has zero infimal Heegaard gradient if and only if 
its fundamental group is infinite. 
\end{theorem}
Remark that 
there are infinitely many Seifert fibered 3-manifolds 
which admit a virtual fibration over 
a surface other than the 2-sphere 
but do not admit that over the circle. 
In fact, it is known that 
such manifolds are precisely the 3-manifolds admitting 
an $\widetilde{SL_2(\mathbb{R})}$-structure. 
See \cite{S} about geometric structures and 
Seifert fibered 3-manifolds for example. \\

To prove Theorem \ref{thm-SF}, 
we will actually construct 
a collection of finite covers with zero Heegaard gradient. 
Our advantage is using some circle bundles over surfaces 
as the finite covers, 
in which Heegaard surfaces are well understood. 
In addition, it will be also shown in Lemma \ref{lem-S1bdl} that 
all minimal genus Heegaard surfaces 
in these covers are \textit{strongly irreducible}. 

The strong irreducibility of Heegaard surface 
was introduced by Casson and Gordon in \cite{CG}. 
The Casson and Gordon's theorem 
is one of the excellent results 
about Haken manifolds and Heegaard splittings: 
if a compact orientable irreducible 3-manifold admits 
an irreducible but weakly reducible Heegaard surface, 
then the manifold is Haken. 
Here a Heegaard surface is called 
\textit{reducible} if it has 
compressing disks on both sides 
whose boundaries coincide. 
A Heegaard surface that is not reducible 
is called \textit{irreducible}. 
A Heegaard surface is called 
\textit{weakly reducible} if it has 
compressing disks on both sides 
which are mutually disjoint. 
A Heegaard surface that is not weakly reducible 
is called \textit{strongly irreducible}. 

Motivated by this result, 
Lackenby also defined in \cite{L} 
the infimal strong Heegaard gradient 
by considering only strongly irreducible Heegaard surfaces 
instead of all Heegaard surfaces. 
Precisely the \textit{infimal strong Heegaard gradient} of 
a collection of finite coverings $\{ M_i \to M \}$ 
of a compact orientable 3-manifold $M$ 
with degree $d_i$ is defined as 
$$ \liminf_i \frac{\chi^{sh}_- (M_i)}{d_i} \ ,$$
where $\chi^{sh}_- (M_i)$ denotes the negative of 
the maximal Euler characteristic of 
a strongly irreducible Heegaard surface in $M_i$. 
Here we set $\chi^{sh}_- (M) = \infty$ 
if $M$ does not have such a Heegaard surface. 
In the same way as before, 
the \textit{infimal strong Heegaard gradient} of $M$ is 
defined as the infimal strong Heegaard gradient of 
the collection of all finite coverings of $M$. 
Again, for brevity, we sometimes drop the word `infimal'. 

Toward the virtually Haken conjecture, 
we are also interested in 
the vanishing of strong Heegaard gradient. 
Concerning this, Lackenby also made a conjecture: 
every closed orientable hyperbolic 3-manifold has 
positive infimal strong Heegaard gradient
(Strong Heegaard gradient conjecture \cite{L}). 
In fact, he proved in \cite{L} that 
a collection of the cyclic covers of 
a closed hyperbolic 3-manifold 
dual to a non-trivial element of $H_2 (M)$ 
has non-zero strong Heegaard gradient. 
In particular, this implies that 
a collection of cyclic covers of 
a closed hyperbolic surface bundle over the circle 
which has zero Heegaard gradient 
has non-zero strong Heegaard gradient. 
On the contrary, 
we show that there are Seifert fibered 3-manifolds 
with zero infimal strong Heegaard gradient. 
\begin{proposition}\label{prop-strong}
A compact orientable Seifert fibered $3$-manifold 
admitting an $\widetilde{SL_2(\mathbb{R})}$-structure 
has zero infimal strong Heegaard gradient. 
\end{proposition}
We also consider the problem 
of which collections of finite coverings 
of a compact orientable Seifert fibered 3-manifold 
have zero Heegaard gradient. 
For a compact orientable hyperbolic 3-manifold $M$ of finite volume, 
Lackenby showed in \cite{L} that 
a collection of the cyclic covers dual to a non-trivial element of 
$H_2 (M, \partial M)$ has zero Heegaard gradient if and only if 
the manifold fibers over the circle and 
the element represents a fiber surface. 
Our result for Seifert fibered 3-manifolds is the following. 
\begin{theorem}\label{thm-cond}
Let $M$ be a closed orientable Seifert fibered $3$-manifold 
with orientable base orbifold of negative Euler characteristic. 
Then the infimal Heegaard gradient of 
a collection $\{ M_i \to M \}$ of finite coverings is zero 
if and only if 
the set of covering degrees between their regular fibers 
induced from $\{ M_i \to M \}$ is unbounded. 
\end{theorem}
Here each cover $M_i$ is endowed with a Seifert fibration 
induced from that of $M$ by covering projection. 
Note that the manifolds considered in Theorem \ref{thm-cond} 
are precisely those admitting 
an $\widetilde{SL_2(\mathbb{R})}$-structure or 
an $\mathbb{H}^2 \times \mathbb{R}$-structure. 
See \cite{S} again. \\

We end this section with an easy observation that 
a compact orientable 3-manifold has negative Heegaard gradient 
if and only if 
it is finitely covered by the 3-sphere $S^3$. 
For let us consider 
the universal covering $\{ S^3 \to M \}$ with degree $d$, 
and then we have $ \chi^h_- ( S^3 ) / d = -2/d < 0 $. 
Together with Theorem \ref{thm-SF}, 
this observation gives the following corollary. 
\begin{corollary}
The Heegaard gradient of a Seifert fibered 3-manifold 
is always non-positive. 
\end{corollary}
\begin{proof}
By Theorem \ref{thm-SF}, Seifert fibered 3-manifolds 
with infinite fundamental group have zero Heegaard gradient. 
On the other hand, it is known that Seifert fibered 3-manifolds 
with finite fundamental group is spherical, 
meaning that finitely covered by the 3-sphere. 
Then they have negative Heegaard gradients as we observed above. 
\end{proof}

\section{Irreducible Heegaard surfaces in 
Seifert fibered $3$-manifolds}\label{sec-SF}

The crucial result we will use throughout this article 
is the classification of 
irreducible Heegaard surfaces of Seifert fibered 3-manifolds. 
It was proved by Moriah and Schultens in \cite{MS} that 
they are either \textit{vertical} or \textit{horizontal}. 

We include here a brief treatment of 
these Heegaard surfaces. 
See \cite{MS,Se} for detailed descriptions. 
In this section, 
let $M$ be a compact orientable Seifert fibered 3-manifold. 
The \textit{base orbifold} $O$ of $M$ is obtained 
by identifying each circle fiber to a point. 
A fiber of $M$ is called 
\textit{regular} or \textit{exceptional} if 
it corresponds to a regular point or a cone point of $O$ 
under the natural projection $M \to O$ respectively. 

First we give 
a construction of a vertical Heegaard surface in $M$. 
For simplicity we suppose that 
$M$ is closed and $O$ is orientable 
with at least two cone points. 
Take a minimal graph $\Gamma$ in the underlying surface $B$ of $O$ 
so that $B - \Gamma$ is a disk and 
$\Gamma$ is disjoint from the cone points of $O$. 
Choose a non-trivial partition of exceptional fibers of $M$. 
For the exceptional fibers in a half of the partition, 
add new edges to $\Gamma$ connecting $\Gamma$ and 
the cone points corresponding to the fibers. 
For all but one of the exceptional fibers in the other half, 
add new loops to $\Gamma$ encircling 
the cone points corresponding to the fibers. 
Lift the graph so obtained to $M$ by the natural projection. 
Add the exceptional fibers in the first half of the partition 
to the lifted graph in $M$, and take a regular neighborhood. 
It can be shown that its boundary surface becomes 
a Heegaard surface of $M$, which is called \textit{vertical}. 
In the case where $M$ has less than two exceptional fibers 
or has non-empty boundary, 
similar construction can apply, and 
we have such a Heegaard surface called vertical. 
Thus all compact orientable Seifert fibered 3-manifold 
with orientable base orbifold contains 
a vertical Heegaard surface. 
Remark that 
the genus of a vertical Heegaard surface 
is independent from the partition chosen. 
It depends only upon 
the number of the exceptional fibers of $M$ 
($=$ the number of the cone points of $O$) and 
the genus of the underlying surface $B$ of $O$. 

Next we consider a horizontal Heegaard surface in $M$. 
Unlike vertical ones, 
it only appears with two restrictions on $M$ as follows. 
Remove the interior of a fibered neighborhood of 
a regular or exceptional fiber from $M$. 
The complement, which remains a Seifert fibered 3-manifold, 
is known to be expressed as a surface bundle over the circle 
with a periodic monodromy. 
The first restriction on $M$ is 
that its fiber surface has a single boundary component. 
Recall that $M$ is reconstructed from the surface bundle 
by gluing the fibered neighborhood, 
which is a solid torus. 
Then the second restriction on $M$ is that 
the meridian of the fibered neighborhood 
intersects the boundary of a fiber surface exactly once. 
Under these restrictions, 
the boundary of the regular neighborhood of 
a fiber surface of the surface bundle yields 
a Heegaard surface in $M$, which is called \textit{horizontal}. 

The result obtained in \cite{MS} makes us possible 
to estimate the \textit{Heegaard genus}, 
i.e., the minimal genus of a Heegaard surface, 
of a closed orientable Seifert fibered 3-manifold as follows. 
\begin{lemma}\label{lem-ineq}
Let $M$ be a closed orientable Seifert fibered 3-manifold 
with orientable base orbifold. 
Then we have 
$$ 2 g(B) + k -2  \leq g (M) \leq 2 g(B) + k +1 ,$$
where $g (M)$ denotes the Heegaard genus of $M$, 
$g(B)$ the genus of the underlying surface $B$ 
of the base orbifold of $M$ and 
$k$ the number of cone points of the base orbifold of $M$. 
\end{lemma}
\begin{proof}
The upper bound is achieved 
by considering a vertical Heegaard surface. 
As stated above, such a Heegaard surface always exists 
and its genus is $2 g(B) + 1$ if $k <2$ and 
$2 g(B) + k -1$ if $k \geq 2$. 
This was originally obtained in \cite{BZ}. 
Typically vertical Heegaard surfaces become minimal genus, 
but in some cases, 
horizontal Heegaard surfaces can be. 
Such cases were studied in \cite{MS}, and 
it is shown in \cite{Se} that 
the Heegaard genus of the manifold 
is one less than that of a vertical Heegaard surface. 
This gives the lower bound of $g(M)$. 
\end{proof}

\section{Vanishing of Heegaard gradient}\label{sec-HG}

In this section, we consider 
the vanishing of Heegaard gradient of Seifert fibered 3-manifolds 
and give a proof of Theorem \ref{thm-SF}. 
The next is the key proposition to prove the theorem. 
\begin{proposition}\label{prop-circle}
An orientable circle bundle over a closed orientable surface 
has zero Heegaard gradient 
if and only if 
its fundamental group is infinite. 
\end{proposition}
\begin{proof}
Let $M$ be an orientable circle bundle over 
a closed orientable surface $F$ of genus $g$ 
and $b(M)$ the obstruction class of $M$, 
which is an integer-valued invariant. 
See \cite{S} for the definition and the properties. 

Assume that $M$ has finite fundamental group. 
Then $M$ is shown to be spherical, 
meaning that $M$ is finitely covered by the 3-sphere. 
Then its Heegaard gradient is negative, 
in particular, is non-zero, as we observed in Section \ref{sec-intro}. 
Remark that 
$F$ is the 2-sphere and $b(M)$ is non-zero in this case. 

Conversely assume that $M$ has infinite fundamental group. 
If $b(M)$ is zero, 
then $M$ is the trivial bundle $F \times S^1$, 
which has zero Heegaard gradient. 
Thus we assume 
that $b(M)$, denoted by $b$ simply, is non-zero. 
By taking the mirror image if necessary, 
we assume that $b >0$. 
Remark that the assumptions that 
$M$ has infinite fundamental group and $b >0$ imply 
that $F$ is not the 2-sphere, i.e., $g \ne 0$. 

Since $g$ is non-zero, 
we can take an $i$-fold covering $F_i \to F$ for $i \geq 2$. 
Note that the Euler characteristic $\chi(F_i)$ of $F_i$ 
is equal to $i (2-2g)$, and so, 
its genus is $1-i(1-g) = ig -i+1$. 
This induces an $i$-fold covering $N_i \to M$, 
where $N_i$ is also a circle bundle over $F_i$. 
Note that the obstruction class of $N_i$ is 
equal to $bi$ \cite[Lemma 3.5]{S}. 

Let $F^0_i$ be the compact surface obtained from $F_i$ 
by removing an open disk. 
Consider the manifold $F^0_i \times S^1$ and 
set the boundary of a surface fiber as a longitude and 
one of the circle fiber as a meridian 
on its boundary torus so that 
$N_i$ is regarded as the 3-manifold 
obtained from $F^0_i \times S^1$ 
by Dehn filling along the slope $bi$. 

Consider the $b i$-fold covering $M_i \to N_i$ 
constructed in the following way. 
First we take the $b i$-fold cyclic covering 
of $F^0_i \times S^1$ in the $S^1$-direction. 
On the boundary of the cover, 
the preimage of the curve of slope $bi$ consists of 
a set of parallel curves of slope 1. 
Let $M_i$ be the Seifert fibered 3-manifold 
obtained from $F^0_i \times S^1$ by Dehn filling along this slope. 
Then it is easily seen that $M_i$ covers $N_i$ 
in the circle direction with degree $b i$. 

The obstruction class of this $M_i$ is equal to $bi/bi = 1$, 
and so, by \cite[Corollary 0.5]{MS}, 
the minimal genus Heegaard surface is horizontal. 
Then the Heegaard genus $g(M_i)$ of $M_i$ is 
equal to twice of the genus of $F_i$, that is, $2(ig -i+1)$. 
By construction, the covering $M_i \to M$ has degree $b i^2$, 
and so 
$$
\inf_i \frac{ \chi^h_- ( M_i ) }{ b i^2 } = 
\inf_i \frac{ 4(ig -i+1) -2 }{ b i^2 } 
\leq \inf_i \frac{ 4 i g  }{ b i^2 } 
= \inf_i \frac{ 4 g }{ b i } = 0 \ .
$$
Since $M$ has infinite fundamental group, 
each $\chi^h_- ( M_i )$ is non-negative, and so, 
the collection $\{ M_i \to M \}$ has zero Heegaard gradient. 
Also, since $M$ has infinite fundamental group, 
the collection of all finite coverings of $M$ 
has non-negative Heegaard gradient, 
and it must be zero as it includes $\{ M_i \to M \}$. 
Thus we conclude that $M$ has zero Heegaard gradient. 
\end{proof}
\begin{lemma}\label{lem-cover}
Let $M$ be a compact orientable $3$-manifold and 
$\widetilde{M}$ a finite cover of $M$. 
If $\widetilde{M}$ has zero Heegaard gradient, 
then also $M$ has. 
\end{lemma}
\begin{proof}
Assume that $\widetilde{M}$ has zero Heegaard gradient. 
Note that $M$, hence $\widetilde{M}$, 
is not finitely covered by the 3-sphere. 
For, if $M$ were, then $M$ and $\widetilde{M}$ 
must have negative Heegaard gradient. 

Let $d$ be the covering degree of $\widetilde{M} \to M$ and 
$\{ \widetilde{M_i} \to \widetilde{M} \}$ 
the collection of all finite covering of $\widetilde{M}$ 
with degree $\widetilde{d_i}$. 
The set of compositions $\{ \widetilde{M_i} \to M \}$ 
gives a collection of finite covering of $M$ 
with degree $\widetilde{d_i} d$. 
Since $M$ is not finitely covered by the 3-sphere, 
any $\chi^h_- ( \widetilde{M_i} )$ is non-negative. 
Thus we have 
$$
0 \leq 
\inf_i \frac{\chi^h_- ( \widetilde{M_i} ) }{ \widetilde{d_i} d }
\leq \inf_i \frac{\chi^h_- ( \widetilde{M_i} ) }{ \widetilde{d_i} } \ . 
$$
Now the right term becomes zero, 
and so $\{ \widetilde{M_i} \to M \}$ has 
zero Heegaard gradient. 
This implies that the Heegaard gradient of 
the set of all finite coverings of $M$ is non-positive. 
But it must be zero as $M$ is not finitely covered by the 3-sphere. 
Therefore $M$ has zero Heegaard gradient. 
\end{proof}
\begin{proof}[Proof of Theorem~{\rm\ref{thm-SF}}]
Let $M$ be a compact orientable Seifert fibered 3-manifold, 
$e(M)$ the Euler number of $M$, 
$O$ the base orbifold of $M$ and 
$\chi(O)$ the Euler characteristic of $O$. 

First suppose that 
$M$ does not virtually fiber over the circle nor 
over a surface other than the 2-sphere. 
Then it is known that $M$ is spherical, 
equivalently 
the fundamental group $\pi_1(M)$ of $M$ is finite. 
As we remarked before, 
it has negative Heegaard gradient, 
in particular, its Heegaard gradient cannot vanish. 
Note that $e(M) \ne 0$ and $\chi(O) > 0$ hold in this case. 

Assume conversely that 
$M$ virtually fibers over the circle or 
over a surface other than the 2-sphere. 
Note that this is equivalent to that 
$M$ has infinite fundamental group. 
In this case, $e(M) = 0$ or $\chi(O) \leq 0$ hold. 

If $M$ virtually fibers over the circle, 
equivalently $e(M)$ is zero, 
then it has zero Heegaard gradient as stated in Section 1. 
Remark that if $M$ has non-empty boundary, 
then $e(M)$ is always zero. 

The remaining case is that 
$M$ virtually fibers over a surface $F$ other than the 2-sphere. 
Since $\chi(O) \leq 0$, 
$F$ has non-positive Euler characteristic, and so, 
a finite cover $\widetilde{M}$ of $M$ fibering over $F$ 
must have infinite fundamental group. 
Taking a double cover if necessary, 
we can assume that $F$ is orientable. 
Then, by Proposition \ref{prop-circle}, 
$\widetilde{M}$ has zero Heegaard gradient. 
This completes the proof of the theorem 
together with Lemma \ref{lem-cover}. 
\end{proof}
In the proof, remark that 
if $e(M) \ne 0$ and $\chi(O) = 0$, 
then $M$ is a Nil-manifold, 
each of which is known to be 
a virtual torus bundle over the circle. 
Precisely, a closed Seifert fibered 3-manifold 
with infinite fundamental group 
does not virtually fiber over the circle 
if and only if 
it admits a geometric structure modeled on 
$\widetilde{SL_2 \mathbb{R} }$. 
See \cite{S} for example. 

\section{$3$-manifolds with zero strong Heegaard gradient}\label{sec-str}

Here we consider an orientable circle bundle 
over a closed orientable surface other than the 2-sphere. 
It is shown to be irreducible and to contain essential tori. 
Thus it is a Haken manifold. 
\begin{lemma}\label{lem-S1bdl}
Let $M$ be an orientable circle bundle over 
a closed orientable surface 
with the obstruction class $\pm 1$. 
Then every irreducible Heegaard surface in $M$ 
is strongly irreducible. 
\end{lemma}
Equivalently such a 3-manifold $M$ contains 
no irreducible but weakly reducible Heegaard surface. 
Thus $M$ is not recognized to be Haken 
by using the Casson and Gordon's theorem. 
\begin{proof}[Proof of Lemma~{\rm\ref{lem-S1bdl}}]
Since the Euler number of $M$ is equal to $\pm 1$, 
$M$ has a unique (up to homeomorphism) irreducible Heegaard surface $S$ 
which is horizontal \cite[Corollary 0.5]{MS}. 
Note that the genus of $S$ is just twice of that of 
the base surface $F$ of $M$. 
Moreover, Corollary 2 in \cite{Sc} says that 
if $S$ is weakly reducible, 
then it is also vertical. 
In our case, a vertical Heegaard surface has genus 
at least twice of that of $F$, and so, 
$S$ must be strongly irreducible. 
\end{proof}
\begin{proof}[Proof of Proposition~{\rm\ref{prop-strong}}]
Let $M$ be a compact orientable Seifert fibered 3-manifold 
which admits an $\widetilde{SL_2(\mathbb{R})}$-structure. 
It is known that $M$ is finitely covered by 
an orientable circle bundle over 
a closed orientable surface 
other than the 2-sphere 
with non-zero obstruction class. 
Then, as in the proofs of 
Theorem \ref{thm-SF} and Lemma \ref{lem-cover}, 
we can find a collection of finite covers 
$\{ M_i \to M \}$ of degree $d_i$ such that 
$$
\inf_i \frac{\chi^h_- ( M_i ) }{ d_i } =0. 
$$
Moreover, by the proof of Proposition \ref{prop-circle}, 
we can take an orientable circle bundle over 
a closed orientable surface 
with obstruction class $\pm 1$ as each $M_i$. 
Then we have $ \chi^h_- ( M_i ) = \chi^{sh}_- ( M_i ) $ holds. 
For Lemma \ref{lem-S1bdl} assures that 
$\chi^h_- ( M_i )$ is attained by 
a strongly irreducible Heegaard surface. 
It implies that the collection of finite covers 
$\{ M_i \to M \}$ has zero strong Heegaard gradient. 
Thus the strong Heegaard gradient of $M$ is non-positive. 
In fact, it must be zero 
because it is at least Heegaard gradient of $M$ 
which is zero as we proved in Theorem \ref{thm-SF}. 
\end{proof}

\section{Vanishing of Heegaard gradient for finite covers} 

In this section, we consider when 
a collection of finite coverings 
of a Seifert fibered 3-manifold has zero Heegaard gradient and 
give a proof of Theorem \ref{thm-cond}. 
\begin{lemma}\label{lem-genus}
Let $M$ be 
a closed orientable irreducible Seifert fibered $3$-manifold 
with infinite fundamental group and orientable base orbifold, and 
$\{ M_i \to M \}$ a collection of finite coverings with degree $d_i$. 
Then the Heegaard gradient of $\{ M_i \to M \}$ is zero 
if and only if 
$ \inf_i g (M_i) / d_i $ is zero, 
where $g (M_i)$ denotes the Heegaard genus of $M_i$. 
\end{lemma}
\begin{proof}
Assume that $ \inf_i g (M_i) / d_i = 0$ holds. 
Then 
$$ 0 \leq \inf_i \frac{\chi^h_- (M_i)}{d_i} 
= \inf_i \frac{ 2 g (M_i) -2 }{d_i}
\leq \inf_i \frac{ 2 g (M_i) }{d_i} =0 \ . $$
Here the first inequality follows from 
the fact that $M$ has infinite fundamental group. 
Thus we have 
the Heegaard gradient of $\{ M_i \to M \}$ is zero. 

Conversely assume that 
the Heegaard gradient of $\{ M_i \to M \}$ is zero. 
We remark that $\chi^h_- (M_i) > 0 $ holds for every $M_i$. 
Otherwise the Heegaard genus $g(M_i)$ is at most 1, 
contradicting the fact that $M$ is irreducible and $\pi_1(M)$ is infinite. 
It follows that $\{ M_i \to M \}$ is an infinite collection and 
the set $\{ d_i \}$ an unbounded set of positive integers. 
Then we have 
\begin{eqnarray*}
0 \leq \inf_i \frac{ g (M_i) }{d_i} 
& \leq & \inf_i \frac{ 2 g (M_i) }{d_i} 
\leq 
\inf_i \left( \frac{ \chi^h_- (M_i) }{d_i} + \frac{2}{d_i} \right). 
\end{eqnarray*}
By the assumption that 
the Heegaard gradient of $\{ M_i \to M \}$ is zero, 
we can find a subsequence $\{ M_j \to M \}$ such that 
$$ \frac{ \chi^h_- (M_j) }{d_j} \to 0 \qquad
\mathrm{as} \quad j \to \infty. $$
Since $\chi^h_- (M_i) > 0 $ holds for every $M_i$, 
we have $ d_j \to \infty$ as $ j \to \infty$. 
Thus we conclude  
$$\left( \frac{ \chi^h_- (M_j) }{d_j} + \frac{2}{d_j} \right) \to 0  
\qquad \mathrm{as} \quad j \to \infty.$$ 
This indicates that 
$$ \inf_i \left( \frac{ \chi^h_- (M_i) }{d_i} + \frac{2}{d_i} \right) 
= 0 ,$$
and the proof is completed. 
\end{proof}
\begin{proof}[Proof of Theorem~{\rm\ref{thm-cond}}]
We first remark that $M$ is irreducible and 
has infinite fundamental group. 
This follows from the assumption that 
the base orbifold of $M$ has negative Euler characteristic. 
See \cite{S} for an example. 

Let $m_i$ be the covering degree between 
the regular fibers of $M_i$ and $M$ 
induced from $\{ M_i \to M \}$. 
Also let $l_i$ be the degree of the orbifold-covering 
between the base orbifolds $O_i$ and $O$ of $M_i$ and $M$ 
induced from $\{ M_i \to M \}$. 
Note that $d_i = m_i l_i$ holds. 

First assume that 
the Heegaard gradient of $\{ M_i \to M \}$ is zero. 
By Lemma \ref{lem-genus}, this implies 
$ \inf_i g (M_i) /d_i = 0$ holds, 
where $g (M_i)$ denotes the Heegaard genus of $M_i$ and 
$d_i$ the covering degree of $M_i \to M$. 
Also by Lemma \ref{lem-ineq}, we have 
$$ g (M_i) \geq 2 g(B_i) + k_i -2 \ ,$$ 
where $g(B_i)$ denotes the genus of 
the underlying surface $B_i$ of $O_i$ and 
$k_i$ the number of cone points of $O_i$. 
Consider the Euler characteristic $\chi (O_i)$ of $O_i$. 
This is calculated as 
$$ \chi (O_i) = 2- 2 g (B_i) - 
\sum_{j=1}^{k_i} \left( 1 - \frac{1}{\alpha_{ij}} \right) \ ,$$ 
where 
$\alpha_{ij}$ denotes the index of the $j$-th cone point of $O_i$. 
Then we have 
$$ 2 g(B_i) + k_i -2 = - \chi(O_i) + \sum_{j=1}^{k_i} \frac{1}{\alpha_{ij}} 
\geq - \chi(O_i) \ .$$
This implies that 
$$ g(M_i) \geq 2 g(B_i) + k_i -2 \geq - \chi(O_i) = - l_i \chi(O) \ ,$$
where $\chi (O)$ denotes the Euler characteristic of $O$. 
Consequently we obtain 
$$  0 = \inf_i \frac{ g (M_i) }{d_i} 
\geq \inf_i \frac{ - l_i \chi(O) }{ m_i l_i} 
= \inf_i \frac{ - \chi(O) }{ m_i } \ .$$ 
By assumption, $\chi(O)$ is negative, and hence, 
the set $\{ m_i \}$ must be unbounded. 

Next we conversely assume that $\{ m_i \}$ is unbounded. 
By Lemma \ref{lem-ineq}, we have
$$ g(M_i) \leq 2 g(B_i) + k_i +1 \ .$$
In the same way as above, 
$$ 2 g(B_i) + k_i +1 = 
- \chi (O_i) + 3 + \sum_{j=1}^{k_i} \frac{1}{\alpha_{ij}} $$ 
holds. 
With $\alpha_{ij} \geq 2$, we obtain
$$ g(M_i) \leq - \chi (O_i) + 3 + \frac{k_i}{2} \ .$$ 
Let $k$ be the number of cone points of $O$. 
Then $k_i$ is at most $l_i k$. 
Thus we conclude 
\begin{eqnarray*}
0 
&& \leq \inf_i \frac{ g (M_i) }{d_i} \leq
\inf_i \frac{ - \chi (O_i) + 3 + k_i/2 }{d_i} \\
&& \leq 
\inf_i \frac{ - l_i \chi (O) + 3 + l_i k/2 }{m_i l_i} 
\leq
\inf_i \left( \frac{ - \chi (O) + k/2 }{m_i} + \frac{3}{m_i l_i} \right). 
\end{eqnarray*}
By the assumption that $\{ m_i \}$ is unbounded, 
we can find a subsequence $\{ M_j \to M \}$ such that 
$m_j \to \infty$ as $j \to \infty$. 
Then, since $l_i \geq 1$ for any $i$, 
$$\left( \frac{ - \chi (O) + k/2 }{m_j} + \frac{3}{m_j l_j} \right) 
\to 0 \qquad \mathrm{as} \quad j \to \infty $$
holds. 
Together with the inequality above, we obtain 
$ \inf_i g (M_i) /d_i = 0$. 
By Lemma \ref{lem-genus}, 
this is equivalent to that 
the Heegaard gradient of $\{ M_i \to M \}$ is zero. 
\end{proof}
%
%
\section*{Acknowledgements}
The author would like to thank 
Yo'av Rieck for helpful discussions. 
He also thanks to Marc Lackenby and Tsuyoshi Kobayashi 
for helpful comments and 
to Eric Sedgwick for suggestions about Proposition \ref{prop-strong}. 
%
%

%
\end{document}